\documentclass{amsart}

\usepackage{qayum_thesis}

\author{Qayum Khan}
\title[On smoothable surgery for 4-manifolds]{On smoothable surgery for 4-manifolds}
\address{Department of Mathematics, Vanderbilt University, Nashville, TN 37240 U.S.A.}
\email{qayum.khan@vanderbilt.edu}

\newcommand{\ncb}[2]{\newcommand{#1}{\mathbf{#2}}}
\ncb{\bK}{K}

\newcommand{\ncs}[2]{\newcommand{#1}{\mathsf{#2}}}
\ncs{\sCW}{CW}

\ncm{\Cl}{Cl} \ncm{\Diffeo}{Diffeo} \ncm{\Homeo}{Homeo}
\ncm{\Hur}{Hur} \ncm{\Torus}{Torus}

\ncq{\qE}{E}

\renewcommand{\CP}{\mathbb{CP}}
\renewcommand{\qE}{\mathbb{E}}
\renewcommand{\F}{\mathbb{F}}
\renewcommand{\qH}{\mathbb{H}}
\renewcommand{\R}{\mathbb{R}}
\renewcommand{\RP}{\mathbb{RP}}
\renewcommand{\Z}{\mathbb{Z}}

\renewcommand{\C}{C}
\newcommand{\CH}{\mathbb{CH}}

\nc{\bbS}{\mathbb{S}}
\ncm{\Out}{Out}

\begin{document}

\begin{abstract}
Under certain homological hypotheses on a compact 4-manifold, we
prove exactness of the topological surgery sequence at the stably
smoothable normal invariants. The main examples are the class of
finite connected sums of 4-manifolds with certain product
geometries. Most of these compact manifolds have non-vanishing
second mod 2 homology and have fundamental groups of exponential
growth, which are not known to be tractable by Freedman-Quinn
topological surgery. Necessarily, the $*$-construction of certain
non-smoothable homotopy equivalences requires surgery on
topologically embedded 2-spheres and is not attacked here by
transversality and cobordism.
\end{abstract}

\maketitle \setcounter{tocdepth}{1} \tableofcontents

\section{Introduction}

\subsection{Objectives}

The main theorem of this paper is a limited form of the surgery
exact sequence for compact 4-manifolds (Theorem \ref{Thm_SES}).
Corollaries include exactness at the smooth normal invariants of the
4-torus $T^4$ (Example \ref{Exm_Torus}) and the real projective
4-space $\RP^4$ (Corollary \ref{Cor_OrderTwo}). C.T.C. Wall proved
an even more limited form of the surgery exact sequence
\cite[Theorem 16.6]{Wall} and remarked that his techniques do not
apply to $T^4$ and $\RP^4$. Although our new hypotheses depend on
the $L$-theory assembly map, we provide a remedy along essentially
the same lines.

\subsection{Results}
Let $(X,\bdry X)$ be a based, compact, connected, topological
4-manifold with fundamental group $\pi = \pi_1(X)$ and orientation
character $\omega = w_1(\tau_X): \pi \to \Z^\x$.  The reader is
referred to Section \ref{Sec_Language} for an explanation of
surgical language.

If $X$ has a preferred smooth structure, consider the following
surgery sequence.
\begin{equation}\label{SES_DIFF_Intro}
\cS^s_\DIFF(X) \xra{\en\eta\en} \cN_\DIFF(X) \xra{\en\sigma_*\en}
L_4^h(\Z[\pi]^\omega)
\end{equation}
Otherwise, let $\TOP 0$ refer to manifolds with the same smoothing
invariant as $X$.
\begin{equation}\label{SES_TOP_Intro}
\cS^s_{\TOP 0}(X) \xra{\en\eta\en} \cN_{\TOP 0}(X)
\xra{\en\sigma_*\en} L_4^h(\Z[\pi]^\omega)
\end{equation}

The first examples consists of orientable 4-manifolds $X$ with
torsion-free, infinite fundamental groups, mostly of exponential
growth.  These include the 4-torus $T^4$ and connected sums of
certain aspherical 4-manifolds of non-positive curvature.

\begin{uncor}[\ref{Cor_FreeProdCrys}]
Let $\pi$ be a free product of groups of the form
\[
\pi = \bigstar_{i=1}^n \Lambda_i
\]
for some $n > 0$, where each $\Lambda_i$ is a torsion-free lattice
in either $\Isom(\qE^{m_i})$ or $\Isom(\qH^{m_i})$ or
$\Isom(\CH^{m_i})$ for some $m_i>0$. Suppose the orientation
character $\omega$ is trivial. Then the surgery sequences
(\ref{SES_DIFF_Intro}) and (\ref{SES_TOP_Intro}) are exact.
\end{uncor}

The second examples consist of a generalization $X$ of
non-aspherical, orientable, simply-connected 4-manifolds. These
include the outcome of smooth surgery on the core circle of the
mapping torus of an orientation-preserving self-diffeomorphism of a
3-dimensional lens space $L(p,q)$. The fundamental groups have
torsion.

\begin{uncor}[\ref{Cor_FreeOdd}]
Let $\pi$ be a free product of groups of the form
\[
\pi = \bigstar_{i=1}^n O_i
\]
for some $n > 0$, where each $O_i$ is an odd-torsion group.
(Necessarily $\omega$ is trivial.) Then the surgery sequences
(\ref{SES_DIFF_Intro}) and (\ref{SES_TOP_Intro}) are exact.
\end{uncor}

The third examples consist of non-aspherical, non-orientable
4-manifolds $X$ whose connected summands are non-orientable with
fundamental group of order two. These include the real projective
4-space $\RP^4$.

\begin{uncor}[\ref{Cor_OrderTwo}]
Suppose $X$ is a $\DIFF$ 4-manifold of the form
\[
X = X_1 \# \cdots \# X_n \# r(S^2 \x S^2)
\]
for some $n > 0$ and $r \geq 0$, and each summand $X_i$ is either
$S^2 \x \RP^2$ or $S^2 \rtimes \RP^2$ or $\#_{S^1} n (\RP^4)$ for
some $1 \leq n \leq 4$. Then the surgery sequences
(\ref{SES_DIFF_Intro}) and (\ref{SES_TOP_Intro}) are exact.
\end{uncor}

The fourth examples consist of orientable 4-manifolds $X$ whose
connected summands are mostly aspherical 3-manifold bundles over the
circle. The important non-aspherical examples include $\# n(S^3\x
S^1)$ with free fundamental group.  The aspherical examples are
composed of fibers of a specific type of Haken 3-manifolds.

\begin{uncor}[\ref{Cor_Haken}]
Suppose $X$ is a $\TOP$ 4-manifold of the form
\[
X = X_1 \# \cdots \# X_n \# r (S^2 \x S^2)
\]
for some $n>0$ and $r \geq 0$, and each summand $X_i$ is the total
space of a fiber bundle
\[
H_i \longra X_i \longra S^1.
\]
Here, we suppose $H_i$ is a compact, connected 3-manifold such that:
\begin{enumerate}
\item
$H_i$ is $S^3$ or $D^3$, or

\item
$H_i$ is irreducible with non-zero first Betti number.
\end{enumerate}
Moreover, if $H_i$ is non-orientable, we assume that the quotient
group $H_1(H_i;\Z)_{(\alpha_i)_*}$ of coinvariants is 2-torsionfree,
where $\alpha_i: H_i \to H_i$ is the monodromy homeomorphism. Then
the surgery sequence (\ref{SES_TOP_Intro}) is exact.
\end{uncor}

Finally, the fifth examples consist of possibly non-orientable
4-manifolds $X$ with torsion-free fundamental group. The connected
summands are surface bundles over surfaces, most of which are
aspherical with fundamental groups of exponential growth. The
aspherical, non-orientable examples of subexponential growth include
simple torus bundles $T^2 \rtimes Kl$ over the Klein bottle,
excluded from Corollary \ref{Cor_FreeProdCrys}.

\begin{uncor}[\ref{Cor_Klein}]
Suppose $X$ is a $\TOP$ 4-manifold of the form
\[
X = X_1 \# \cdots \# X_n \# r(S^2 \x S^2)
\]
for some $n > 0$ and $r \geq 0$, and each summand $X_i$ is the total
space of a fiber bundle
\[
\Sigma_i^f \longra X_i \longra \Sigma_i^b.
\]
Here, we suppose the fiber and base are compact, connected
2-manifolds, $\Sigma_i^f \neq \RP^2$, and $\Sigma_i^b$ has positive
genus. Moreover, if $X_i$ is non-orientable, we assume that the
fiber $\Sigma_i^f$ is orientable and that the monodromy action of
$\pi_1(\Sigma_i^b)$ of the base preserves any orientation on the
fiber (i.e. the bundle is simple). Then the surgery sequence
(\ref{SES_TOP_Intro}) is exact.
\end{uncor}

\subsection{Techniques}

Our methods employ various bits of geometric topology: topological
transversality in all dimensions (Freedman-Quinn \cite{FQ}), and the
analysis of smooth normal invariants of the Novikov pinching trick,
which is used to construct homotopy self-equivalences of 4-manifolds
(Cochran-Habegger, Wall \cite{CH}). Our hypotheses are
algebraic-topological in nature and come from the surgery
characteristic class formulas of Sullivan-Wall \cite{Wall} and from
the assembly map components of Taylor-Williams \cite{TW}, as well as
control of $\pi_2$ in non-orientable cases.

Jonathan Hillman has successfully employed these now standard
techniques to classify 4-manifolds, up to $s$-cobordism, in the
homotopy type of certain surface bundles over surfaces \cite[\S
2]{Hillman_SurfaceBundles} \cite[Ch. 6]{Hillman4Book}. Along the
same lines, our abundant families of 4-manifold examples also have
fundamental groups of exponential growth, and so, too, are currently
inaccessible by topological surgery \cite{FQ, FreedmanTeichner,
KrushkalQuinn}.

The reader should be aware that the topological transversality used
in Section \ref{Sec_Surgery4} produces 5-dimensional $\TOP$ normal
bordisms $W\to X \x \Delta^1$ which may not be smoothable, although
the boundary $\bdry W = \bdry_- W \cup \bdry_+ W$ is smoothable. In
particular, $W$ may not admit a $\TOP$ handlebody structure relative
to $\bdry_- W$. Hence $W$ may not be the trace of surgeries on
topologically embedded 2-spheres in $X$. Therefore, in general, $W$
cannot be produced by Freedman-Quinn surgery theory, which has been
developed only for a certain class of fundamental groups $\pi_1(X)$
of subexponential growth. In this way, topological cobordism is
superior to surgery.

\subsection{Language}\label{Sec_Language}

For any group $\pi$, we shall write $\Wh_0(\pi) :=
\wt{K}_0(\Z[\pi])$ for the \textbf{projective class group} and
$\Wh_1(\pi) := \wt{K}_1(\Z[\pi])/\gens{\pi}$ for the
\textbf{Whitehead group}.

Let $\CAT$ be either the manifold category $\TOP$ or $\PL=\DIFF$ in
dimensions $< 7$. Suppose $(X,\bdry X)$ is a based, compact,
connected $\CAT$ 4-manifold. Let us briefly introduce some basic
notation used throughout this paper. The fundamental group $\pi =
\pi_1(X)$ depends on a choice of basepoint; a basepoint is essential
if $X$ is non-orientable.  The \textbf{orientation character}
$\omega = w_1(X): \pi \to \Z^\x$ is a homomorphism that assigns $+1$
or $-1$ to a loop $\lambda: S^1 \to X$ if the pullback bundle
$\lambda^*(\tau_X)$ is orientable or non-orientable. Recall that any
finitely presented group $\pi$ and arbitrary orientation character
$\omega$ can be realized on some closed, smooth 4-manifold $X$ by a
straightforward surgical construction. A choice of generator $[X]
\in H_4(X,\bdry X;\Z^\omega)$ is called a \textbf{twisted
orientation class}.

Let us introduce the terms in the surgery sequence investigated in
Section \ref{Sec_Surgery4}. The \textbf{simple structure set}
$\cS^s_\CAT(X)$ consists of $\CAT$ $s$-bordism classes in
$\R^\infty$ of simple homotopy equivalences $h: Y \to X$ such that
$\bdry h: \bdry Y \to \bdry X$ is the identity. Here,
\textbf{simple} means that the torsion of the acyclic
$\Z[\pi]$-module complex $\mathrm{Cone}(\wt{h})$ is zero, for some
preferred finite homotopy CW-structures on $Y$ and $X$. Indeed, any
compact topological manifold $(X,\bdry X)$ has a canonical simple
homotopy type, obtained by cleanly embedding $X$ into euclidean
space \cite[Thm. III.4.1]{KS}. Therefore, for any abelian group $A$
and $n>1$, pulling back the inverse of the Hurewicz isomorphism
induces a bijection from $[X/\bdry X, K(A,n)]_0$ to $H^n(X,\bdry
X;A)$. This identification shall be used implicitly throughout the
paper.

Denote $\G_n$ as the topological monoid of homotopy
self-equivalences $S^{n-1} \to S^{n-1}$, and $\G := \colim_n \G_n$
as the direct limit of $\set{ \G_n \to \G_{n+1}}$. The
\textbf{normal invariant set} $\cN_\CAT(X) \iso [X/\bdry X,
\G/\CAT]_0$ consists of $\CAT$ normal bordism classes in $\R^\infty$
of degree one, $\CAT$ normal maps $f: M \to X$ such that $\bdry f:
\bdry X \to \bdry X$ is the identity; we suppress the normal data
and define $X/\varnothing = X \sqcup \pt$. Denote $\hat{f}: X/\bdry
X \to \G/\CAT$ as the associated homotopy class of based maps.
Indeed, transversality in the $\TOP$ category holds for all
dimensions and codimensions \cite{KS, FQ}. The normal invariants map
$\eta: \cS^s_\CAT(X) \to \cN_\CAT(X)$ is a forgetful map. The
\textbf{surgery obstruction group} $L_4^h(\Z[\pi]^\omega)$ consists
of Witt classes of nonsingular quadratic forms over the group ring
$\Z[\pi]$ with involution $(g \mapsto \omega(g) g\inv)$. The surgery
obstruction map $\sigma_*^h: \cN_\CAT(X) \to L_4^h(\Z[\pi]^\omega)$
vanishes on the image of $\eta$. The basepoint of the former two
sets is the identity map $\Id_X: X \to X$, and the basepoint of the
latter set is the Witt class $0$.

\subsection{Invariants}\label{Subsec_Preliminaries}

The unique homotopy class of classifying maps $u: X \to B\pi$ of the
universal cover induces homomorphisms
\begin{gather*}
u_0: H_0(X;\Z^\omega) \longra H_0(\pi;\Z^\omega)\\
u_2: H_2(X;\Z_2) \longra H_2(\pi;\Z_2).
\end{gather*}

Next, recall that the manifold $X$ has a second Wu class
\[
v_2(X) \in H^2(X;\Z_2) = \Hom(H_2(X;\Z_2),\Z_2)
\]
defined for all $a \in H^2(X,\bdry X;\Z_2)$ by $\gens{v_2(X), a \cap
[X]} = \gens{a \cup a, [X]}$. This unoriented cobordism
characteristic class is uniquely determined from the Stiefel-Whitney
classes of the tangent microbundle $\tau_X$ by the formula
\[
v_2(X) = w_1(X) \cup w_1(X) + w_2(X).
\]
Observe that $v_2(X)$ vanishes if $X$ is a $\TOP$ Spin-manifold.

Finally, let us introduce the relevant surgery characteristic
classes. Observe
\begin{equation}\label{Eqn_H0}
H_0(\pi;\Z^\omega) =
\Z/\gens{\omega(g) - 1 \ST g \in \pi} =
\begin{cases}\Z & \text{if } \omega = 1\\ \Z_2 & \text{if } \omega \neq 1.\end{cases}
\end{equation}
The 0th component of the 2-local assembly map $A_\pi$ \cite{TW} has
an integral lift
\[
I_0 : H_0(\pi;\Z^\omega) \longra L_4^h(\Z[\pi]^\omega).
\]
The image $I_0(1)$ equals the Witt class of the $E_8$ quadratic form
\cite[Rmk. 3.7]{Davis}. The 2nd component of the 2-local assembly
map $A_\pi$ \cite{TW} has an integral lift
\[
\kappa_2: H_2(\pi;\Z_2) \longra L_4^h(\Z[\pi]^\omega).
\]

Let $f: M \to X$ be a degree one, $\TOP$ normal map. According to
Ren\'e Thom \cite{Thom}, every homology class in $H_2(X,\bdry
X;\Z_2)$ is represented by $g_*[\Sigma]$ for some compact, possibly
non-orientable surface $\Sigma$ and $\TOP$ immersion $g:
(\Sigma,\bdry \Sigma) \to (X,\bdry X)$.  The codimension two
Kervaire-Arf invariant
\[
\kerv(f): H_2(X,\bdry X;\Z_2) \longra \Z_2
\]
assigns to each two-dimensional homology class $g_*[\Sigma]$ the Arf
invariant of the degree one, normal map $g^*(f): f^*(\Sigma) \to
\Sigma$. The element $\kerv(f) \in H^2(X,\bdry X;\Z_2)$ is invariant
under $\TOP$ normal bordism of $f$; it may not vanish for homotopy
equivalences. If $M$ and $X$ are oriented, then there is a signature
invariant
\[
\sign(f) := (\sign(M) - \sign(X))/8 \in H_0(X;\Z),
\]
which does vanish for homotopy equivalences.  For any compact
topological manifold $X$, the Kirby-Siebenmann invariant $\ks(X) \in
H^4(X,\bdry X;\Z_2)$ is the sole obstruction to the existence of a
$\DIFF$ structure on $X \x \R$ or equivalently on $X \# r(S^2 \x
S^2)$ for some $r \geq 0$. Furthermore, the image of $\ks(X) \cap
[X]$ in $\Z_2$ under the augmentation map $X \to \pt$ is an
invariant of unoriented $\TOP$ cobordism \cite[\S10.2B]{FQ}. Define
\[
\ks(f) := f_*(\ks(M) \cap [M]) - (\ks(X) \cap [X]) \in H_0(X;\Z_2).
\]

In Section \ref{Sec_Surgery4}, we shall use Sullivan's surgery
characteristic class formulas as geometrically identified in
dimension four by J.F. Davis \cite[Prop. 3.6]{Davis}:
\begin{eqnarray}\label{Eqn_Invariants1}
\hat{f}^*(k_2) \cap [X] &=& \kerv(f) \cap [X] \in H_2(X;\Z_2)\\
\label{Eqn_Invariants2} \hat{f}^*(\ell_4) \cap [X] &=&
\begin{cases}\sign(f) \in H_0(X;\Z)
& \text{if } \omega = 1\\
\ks(f) + (\kerv(f)^2 \cap [X]) \in H_0(X;\Z_2) & \text{if } \omega
\neq 1.
\end{cases}
\end{eqnarray}
Herein is used the 5th stage Postnikov tower \cite{KS, Wall}
\[
k_2 + \ell_4: \G/\TOP^{[5]} \xra{\en\simeq\en} K(\Z_2,2) \x K(\Z,4).
\]
The two expressions in (\ref{Eqn_Invariants2}) agree modulo two
\cite[Annex 3, Thm. 15.1]{KS}:
\begin{equation}\label{Eqn_Invariants3}
\ks(f) = \prn{\red_2 (\hat{f})^*(\ell_4) - (\hat{f})^*(k_2)^2} \cap
[X] \in H_0(X;\Z_2).
\end{equation}

\section{Smoothing normal bordisms}

Let $(X,\bdry X)$ be a based, compact, connected, $\DIFF$
4-manifold. We start with group-theoretic criteria on the existence
and uniqueness of smoothing the topological normal bordisms relative
$\bdry X$ from the identity map on $X$ to itself.

\begin{prop}\label{Prop_RedTOP}
With respect to the Whitney sum $H$-space structures on the $\CAT$
normal invariants, there are exact sequences of abelian groups:
\begin{center}\tiny
$0 \longra \Tor_1( H_0(\pi;\Z^\omega), \Z_2 ) \longra \cN_\DIFF(X)
\xra{\red_\TOP} \cN_\TOP(X) \xra{\en\ks\en} H_0(\pi;\Z_2) \longra 0
$ \end{center}
\begin{center}\tiny
$0 \longra \Tor_1( H_1(\pi;\Z^\omega), \Z_2) \longra \cN_\DIFF(X \x
\Delta^1) \xra{\red_\TOP} \cN_\TOP(X \x \Delta^1) \xra{\en\ks\en}
H_1(\pi;\Z^\omega) \xo \Z_2 \longra 0.$
\end{center}
\end{prop}

\begin{proof}
Since $X$ is a $\CAT$ manifold, by $\CAT$ transversality and Cerf's
result that $\PL/\mathrm{O}$ is 6-connected \cite{KS,FQ}, we can
identify the based sets
\begin{eqnarray*}
\cN_\DIFF(X) &=& [X/\bdry X, \G/\PL]_0\\
\cN_\TOP(X) &=& [X/\bdry X, \G/\TOP]_0\\
\cN_\DIFF(X \x \Delta^1) &=& [S^1 \wedge (X/\bdry X), \G/\PL]_0\\
\cN_\TOP(X \x \Delta^1) &=& [S^1 \wedge (X/\bdry X), \G/\TOP]_0.
\end{eqnarray*}
Furthermore, each right-hand set is an abelian group with respect to
the $H$-space structure on $\G/\CAT$ given by Whitney sum of $\CAT$
microbundles.

For any based space $Z$ with the homotopy type of a CW-complex,
there is the Siebenmann-Morita exact sequence of abelian groups
\cite[Annex 3, Thm. 15.1]{KS}:
\begin{center}\tiny
$0 \longra \Cok\prn{\red_2^{(3)}} \longra  [Z,\G/\PL]_0
\xra{\red_\TOP} [Z, \G/\TOP]_0  \xra{\en\ks\en}
\Img\prn{\red_2^{(4)} + \Sq^2} \longra 0.$
\end{center}
Here, the stable cohomology operations
\begin{gather*}
\red_2^{(n)}: H^n(Z;\Z) \longra H^n(Z;\Z_2) \\
\Sq^2: H^2(Z;\Z_2) \to H^4(Z;\Z_2)
\end{gather*}
are reduction modulo two and the second Steenrod square. The
homomorphism $\ks$ is given by the formula $\ks(a,b) =
\red_2^{(4)}(a) - \Sq^2(b)$, as stated in (\ref{Eqn_Invariants3}),
which follows from Sullivan's determination (\ref{Eqn_Sullivan})
below.

Suppose $Z = X/\bdry X$. By Poincar\'e duality and the universal
coefficient sequence, there are isomorphisms
\begin{gather*}
\Cok\prn{\red_2^{(3)}} \iso \Cok\prn{\red_2: H_1(X;\Z^\omega) \to
H_1(X;\Z_2)} \iso \Tor_1(H_0(X;\Z^\omega), \Z_2)\\
\Img\prn{\red_2^{(4)}} \iso \Img\prn{\red_2: H_0(X;\Z^\omega) \to
H_0(X;\Z_2)} = H_0(X;\Z_2).
\end{gather*}
Therefore we obtain the exact sequence for the normal invariants of
$X$.

Suppose $Z = S^1 \wedge (X/\bdry X)$. By the suspension isomorphism
$\Sigma$, Poincar\'e duality, and the universal coefficient
sequence, there are isomorphisms
\begin{multline*}
\Cok\prn{\red_2^{(3)}} \iso \Cok\prn{\red_2^{(2)}: H^2(X,\bdry X;\Z)
\to H^2(X,\bdry X; \Z_2)}\\ \iso \Cok\prn{\red_2: H_2(X;\Z^\omega)
\to H_2(X;\Z_2)} \iso \Tor_1(H_1(X;\Z^\omega), \Z_2).
\end{multline*}
Note, since the cohomology operations $\red_2^{(4)}$ and $\Sq^2$ are
stable, that
\[
(\Sigma\inv \circ \ks)(\Sigma a, \Sigma b) = \red_2^{(3)}(a) -
\Sq^2(b) = \red_2^{(3)}(a)
\]
for all $a \in H^3(X,\bdry X;\Z)$ and $b \in H^1(X,\bdry X;\Z_2)$.
Then, by Poincar\'e duality and the universal coefficient sequence,
we have
\begin{multline*}
\Img(\ks) \iso \Img\prn{\red_2^{(3)}: H^3(X,\bdry X;\Z) \to
H^3(X,\bdry X; \Z_2)}\\
\iso \Img\prn{\red_2: H_1(X;\Z^\omega) \to H_1(X;\Z_2)} =
H_1(\pi;\Z^\omega) \xo \Z_2.
\end{multline*}
Therefore we obtain the exact sequence for the normal invariants of
$X \x \Delta^1$.
\end{proof}

\section{Homotopy self-equivalences}\label{Sec_SHE}

Recall Sullivan's determination \cite{MM, Wall}
\begin{equation}\label{Eqn_Sullivan}
k_2 + 2 \ell_4: \G/\PL^{[5]} \xra{\en\simeq\en} K(\Z_2,2)
\x_{\delta(\Sq^2)} K(\Z,4).
\end{equation}
The homomorphism $\delta: H^4(K(\Z_2,2);\Z_2) \to H^5(K(\Z_2,2);\Z)$
is the Bockstein associated to the coefficient exact sequence $0 \to
\Z \xra{2} \Z \to \Z_2 \to 0$, and the element $\Sq^2$ is the 2nd
Steenrod square.  The cohomology classes $k_2$ and $2 \ell_4$ map to
a generator of the base $K(\Z_2,2)$ and of the fiber $K(\Z,4)$.
Moreover, the cohomology class $2 \ell_4$ of $\G/\PL$ is the
pullback of the cohomology class $\ell_4$ of $\G/\TOP$ under the
forgetful map $\red_\TOP: \G/\PL \to \G/\TOP$.  The above homotopy
equivalence gives
\[
\red_2(2\ell_4) = (k_2)^2 \in H^4(\G/\PL;\Z_2);
\]
compare \cite[Theorem 4.32, Footnote]{MM}. There exists a symmetric
$L$-theory twisted orientation class $[X]_{\qL^\cdot} \in
H_4(X,\bdry X; \qL^{\cdot\omega})$ fitting into a commutative
diagram
\begin{center}\small
$\xymatrix{ \cN_\PL(X) \ar[rr]^{\sigma_*} \ar@{=}[d] & &
L_4^h(\Z[\pi]^\omega)\\ [X/ \bdry X, \G/\PL]_0 \ar[r]
\ar[d]_{\red_\TOP} & \wt{\Omega}_4^\SPL( B\pi_+\wedge \G/\PL,\omega)
\ar[ur] \ar[d]_{\red_\TOP} & H_0(\pi; \Z^\omega) \oplus H_2(\pi;
\Z_2) \ar[u]_{A_\pi\!\gens{1}} \ar@{=}[d] \\ [X/\bdry X,\G/\TOP]_0
\ar@{=}[d] \ar[r] & \wt{H}_4(B\pi_+\wedge \G/\TOP; \MSTOP^\omega)
\ar[r]^{\qquad\mathrm{act}_*} & H_4(B\pi; \GTOP^\omega )\\ H^0(X,
\bdry X; \GTOP) \ar[r]^{\cap [X]_{\qL^\cdot}} & H_4(X;\GTOP^\omega)
\ar[ru]_{u_*} &  }$
\end{center}\vspace{2mm}
due to Sullivan-Wall \cite[Thm. 13B.3]{Wall} and Quinn-Ranicki
\cite[Thm. 18.5]{RanickiTOP}. Here, the identification $\cN_\PL(X) =
[X/\bdry X,\G/\PL]_0$ only makes sense if $\ks(X)=0$. It follows
that the image $\hat{\sigma}(g) \in H_4(\pi; \GTOP^\omega)$, through
the scalar product $\mathrm{act}$\footnote{$\qL.\gens{1} = \GTOP$ is
a module spectrum over the ring spectrum $\qL^\cdot = \MSTOP$ via
Brown representation. Refer to \cite[Rmk. B9]{RanickiTOP} \cite[Thm.
9.8]{Wall} on the level of homotopy groups or to Sullivan's original
method of proof in his thesis \cite{Sullivan}.}, of a normal
invariant $g: X / \bdry X \to \G/\PL$ consists of two characteristic
classes:
\[
\hat{\sigma}(g) = u_0\prn{g^*(2\ell_4)\cap[X]} \oplus
u_2\prn{g^*(k_2)\cap[X]},
\]
which are determined by the $\TOP$ manifold-theoretic invariants in
Subsection \ref{Subsec_Preliminaries}.

We caution the reader that $\ell_4 \notin H^4(\G/\PL;\Z)$; the
notation $2\ell_4$ is purely formal.

\begin{defn}
Let $(X,\bdry X)$ be any based, compact, connected $\TOP$
$4$-manifold. Define the \textbf{stably smoothable} subsets
\begin{eqnarray*}
\cN_{\TOP 0}(X) &:=& \set{f \in \cN_\TOP(X) \ST \ks(f)=0 }\\
\cS^s_{\TOP 0}(X) &:=& \set{h \in \cS^s_\TOP(X) \ST \ks(h)=0 }.
\end{eqnarray*}
\end{defn}

Recall $X$ has fundamental group $\pi$ and orientation character
$\omega$.

\begin{hyp}\label{Hyp_Orient} Let $X$ be orientable.
Suppose that the homomorphism
\[
\kappa_2: H_2(\pi;\Z_2) \longra L_4^h(\Z[\pi]^\omega)
\]
is injective on the subgroup $u_2(\Ker\,v_2(X))$.
\end{hyp}

\begin{hyp}\label{Hyp_NonorientFin}
Let $X$ be non-orientable such that $\pi$ contains an
orientation-reversing element of finite order, and if $\CAT=\DIFF$,
then suppose that orientation-reversing element has order two.
Suppose that $\kappa_2$ is injective on all $H_2(\pi;\Z_2)$, and
suppose that $\Ker(u_2) \subseteq \Ker(v_2)$.
\end{hyp}

\begin{hyp}\label{Hyp_NonorientInf}
Let $X$ be non-orientable such that there exists an epimorphism
$\pi^\omega \to \Z^-$. Suppose that $\kappa_2$ is injective on the
subgroup $u_2(\Ker\,v_2(X))$.
\end{hyp}

\begin{prop}\label{Prop_SelfEquivalence}
Let $f: M \to X$ be a degree one, normal map of compact, connected
$\TOP$ 4-manifolds such that $\bdry f = \Id_{\bdry X}$. Suppose
Hypothesis \ref{Hyp_Orient} or \ref{Hyp_NonorientFin} or
\ref{Hyp_NonorientInf}. If $\sigma_*(f)=0$ and $\ks(f)=0$, then $f$
is $\TOP$ normally bordant to a homotopy self-equivalence $h: X \to
X$ relative to $\bdry X$.
\end{prop}

\begin{proof}
Since $\ks(f)=0$, there is a (formal) based map $g: X/\bdry X \to
\G/\PL$ such that
\[
\red_\TOP \circ g = \hat{f} : X/\bdry X \to \G/\TOP.
\]
So $g$ has vanishing surgery obstruction:
\[
0 = \sigma_*(f) = (I_0 + \kappa_2) (\hat{\sigma}(g)) \in
L_4^h(\Z[\pi]^\omega).
\]

\emph{Suppose $X$ is orientable; that is, $\omega=1$.}  Then the
inclusion $1^+ \to \pi^\omega$ is retractive and induces a split
monomorphism $L_4^h(\Z[1]) \to L_4^h(\Z[\pi])$ with cokernel defined
as $\tilde{L}_4^h(\Z[\pi])$. So the above sum of maps is direct:
\[
0 = (I_0 \oplus \kappa_2)(\hat{\sigma}(g)) \in L_4^h(\Z[1]) \oplus
\tilde{L}_4^h(\Z[\pi]).
\]
Then both the signature and the square of the Kervaire-Arf invariant
vanish (\ref{Eqn_Invariants3}):
\[
0 = g^*(2\ell_4) = g^*(k_2)^2.
\]
So $(\hat{f})^*(k_2) \cap [X] \in \Ker\, v_2(X)$. Therefore, since
$\kappa_2$ is injective on the subgroup $u_2(\Ker\, v_2(X))$, we
have $(\hat{f})^*(k_2) \cap [X] \in \Ker(u_2)$. So, by the Hopf
exact sequence
\[
\pi_2(X) \xo \Z_2 \xra{\Hur} H_2(X;\Z_2) \xra{\en u_2 \en}
H_2(\pi;\Z_2) \longra 0,
\]
there exists $\alpha \in \pi_2(X)$ such that our homology class is
spherical:
\[
(\hat{f})^*(k_2) \cap [X] = (\red_2 \circ \Hur)(\alpha).
\]

\emph{Suppose $X$ is non-orientable; that is, $\omega \neq 1$.} Let
$x \in \pi$ be an orientation-reversing element: $\omega(x) = -1$.
First, consider the case that $x$ has finite order. By taking an odd
order power, we may assume that $x$ has order $2^N$ for some $N
> 0$.  Then the map induced by $1^+ \to
\pi^\omega$ has a factorization through $(\C_{2^N})^-$:
\[
L_4^h(\Z[1]) \longra L_4^h(\Z[\C_{2^N}]^-) \xra{\en x_* \en}
L_4^h(\Z[\pi]^\omega).
\]
The abelian group in the middle is zero by \cite[Theorem 3.4.5,
Remark]{WallGroupRing}.  Then, since $I_0: H_0(\pi;\Z^\omega) \to
L_4^h(\Z[\pi]^\omega)$ factors through $L_4^h(\Z[1]) \iso \Z$,
generated by the Witt class $[E_8]$, we must have $I_0 = 0$. So
\[
0 = \sigma_*(f) = \kappa_2 (\hat{\sigma}(f)) \in
L_4^h(\Z[\pi]^\omega),
\]
and since $\kappa_2$ is injective on all $H_2(\pi;\Z_2)$, we have
\[
0 = \hat{\sigma}(f) = u_2\prn{(\hat{f})^*(k_2) \cap [X]}.
\]
Then $(\hat{f})^*(k_2) \cap [X] \in \Ker(u_2) \subseteq \Ker(v_2)$
by hypothesis, and the class is spherical.

Next, consider the case there are no orientation-reversing elements
of finite order.  Then, by hypothesis, there is an epimorphism $p:
\pi^\omega \to \Z^-$, which is split by a monomorphism with image
generated by some orientation-reversing infinite cyclic element $y
\in \pi$.  Define $\ol{L}_4^h(\Z[\pi]^\omega)$ as the kernel of
$p_*$. Then $y_*$ induces a direct sum decomposition
\[
L_4^h(\Z[\pi]^\omega) = L_4^h(\Z[\Z]^-) \oplus
\ol{L}_4^h(\Z[\pi]^\omega).
\]
The abelian group of the non-orientable Laurent extension in the
middle is isomorphic to $\Z_2$, generated by the Witt class $[E_8]$,
according to the quadratic version of \cite[Theorem 4.1]{MR2} with
orientation $u=-1$.  Then the map $I_0: H_0(\pi;\Z^\omega) \to
L_4^h(\Z[\pi]^\omega)$ factors through the summand $L_4^h(\Z[\Z]^-)$
by an isomorphism; functorially, $\kappa_2$ has zero projection onto
that factor. So the sum of maps is direct, similar to the oriented
case:
\[
0 = (I_0 \oplus \kappa_2)(\hat{\sigma}(g)) \in L_4^h(\Z[\Z]^-)
\oplus \ol{L}_4^h(\Z[\pi]^\omega).
\]
A similar argument, using the smooth normal invariant $g$, shows
that
\[
0 = (\hat{f})^*(\ell_4) =  (\hat{f})^*(k_2)^2.
\]
Hence $(\hat{f})^*(k_2) \cap [X] \in \Ker\,v_2(X)$.  Since
$\kappa_2$ is injective on $u_2(\Ker\,v_2(X))$, we also have
$(\hat{f})^*(k_2) \cap [X] \in \Ker(u_2)$, thus the class is
spherical.

\emph{Let us return to the general case of $X$ without any condition
on orientability.} For any $\alpha \in\pi_2(X)$, there is a homotopy
operation, called the Novikov pinch map, defined by the homotopy
self-equivalence
\[
h: X \xra{\pinch} X \vee S^4 \xra{\Id_X \vee \Sigma\eta} X \vee S^3
\xra{\Id_X \vee \eta} X \vee S^2 \xra{\Id_X \vee \alpha} X.
\]
Here, $\eta: S^3 \to S^2$ and $\Sigma \eta: S^4 \to S^3$ are the
complex Hopf map and its suspension that generate the stable
homotopy groups $\pi_1^s$ and $\pi_2^s$.

For the normal invariant of the self-equivalence $h: X \to X$
associated to our particular $\alpha$, there is a formula in the
simply-connected case due to Cochran and Habegger \cite[Thm.
5.1]{CH} and generalized to the non-simply connected case by Kirby
and Taylor \cite[Thm. 18, Remarks]{KT}:
\begin{gather*}
(\hat{h})^*(k_2) = \prn{1 + \gens{v_2(X), (\hat{f})^*(k_2) \cap
[X]}} \cdot (\hat{f})^*(k_2) = (\hat{f})^*(k_2)\\
(\hat{h})^*(\ell_4) = 0 = (\hat{f})^*(\ell_4).
\end{gather*}
Here, we have used $(\hat{f})^*(k_2) \cap [X] \in \Ker\,v_2(X)$ and,
if $X$ is non-orientable, $\ks(f)=0$ in (\ref{Eqn_Invariants2}).
Therefore $f: M \to X$ is $\TOP$ normally bordant to the homotopy
self-equivalence $h: X \to X$ relative to the identity $\bdry X \to
\bdry X$ on the boundary.
\end{proof}

\section{Smoothable surgery for 4-manifolds}\label{Sec_Surgery4}

Terry Wall asked if the smooth surgery sequence is exact at the
normal invariants for the 4-torus $T^4$ and real projective 4-space
$\RP^4$; see the remark after \cite[Thm. 16.6]{Wall}. The latter
case of $\RP^4$ was affirmed implicitly in the work of Cappell and
Shaneson \cite{CS2}. The main theorem of this section affirms the
former case of $T^4$ and extends their circle sum technique for
$\RP^4$ to a broader class of non-orientable 4-manifolds, using the
assembly map and smoothing theory.

\begin{thm}\label{Thm_SES}
Let $(X,\bdry X)$ be a based, compact, connected, $\CAT$ 4-manifold
with fundamental group $\pi = \pi_1(X)$ and orientation character
$\omega = w_1(X): \pi \to \Z^\x$.
\begin{enumerate}
\item
Suppose Hypothesis \ref{Hyp_Orient} or \ref{Hyp_NonorientFin}. Then
the surgery sequence of based sets is exact at the smooth normal
invariants:
\begin{equation}\label{SES_DIFF}
\cS^s_\DIFF(X) \xra{\en\eta\en} \cN_\DIFF(X) \xra{\en\sigma_*\en}
L_4^h(\Z[\pi]^\omega).
\end{equation}
\item
Suppose Hypothesis \ref{Hyp_Orient} or \ref{Hyp_NonorientFin} or
\ref{Hyp_NonorientInf}.  Then the surgery sequence of based sets is
exact at the stably smoothable normal invariants:
\begin{equation}\label{SES_TOP}
\cS^s_{\TOP 0}(X) \xra{\en\eta\en} \cN_{\TOP 0}(X)
\xra{\en\sigma_*\en} L_4^h(\Z[\pi]^\omega).
\end{equation}
\end{enumerate}
\end{thm}

The above theorem generalizes a statement of Wall \cite[Theorem
16.6]{Wall} proven correctly by Cochran and Habegger \cite{CH} for
closed, oriented $\DIFF$ 4-manifolds.

\begin{cor}[Wall]\label{Cor_Wall}
Suppose the orientation character $\omega$ is trivial and the group
homology vanishes: $H_2(\pi;\Z_2) = 0$.  Then the surgery sequence
(\ref{SES_DIFF}) is exact.
\end{cor}

A fundamental result from geometric group theory is that any
torsion-free, finitely presented group $\Gamma$ is of the form
$\Gamma = \bigstar_{i=1}^n \Gamma_i$ for some $n \geq 0$, where each
$\Gamma_i$ is either $\Z$ or a one-ended, finitely presented group.
Geometric examples of such $\Gamma_i$ are torsion-free lattices of
any rank. The Borel/Novikov Conjecture (i.e. Integral Novikov
Conjecture) would imply that $\kappa_2$ is injective for all
finitely generated, torsion-free groups $\pi$ and all $\omega$
\cite{Davis}.  At the moment, we have:

\begin{cor}\label{Cor_FreeProdCrys}
Let $\pi$ be a free product of groups of the form
\[
\pi = \bigstar_{i=1}^n \Lambda_i
\]
for some $n > 0$, where each $\Lambda_i$ is a torsion-free lattice
in either $\Isom(\qE^{m_i})$ or $\Isom(\qH^{m_i})$ or
$\Isom(\CH^{m_i})$ for some $m_i>0$. Suppose the orientation
character $\omega$ is trivial. Then the surgery sequences
(\ref{SES_DIFF}) and (\ref{SES_TOP}) are exact.
\end{cor}

\begin{exm}\label{Exm_Torus}
Besides stabilization with connected summands of $S^2 \x S^2$, the
preceding corollary includes the orientable manifolds $X = T^4 =
\prod 4 (S^1)$ and $X = \# n (S^1 \x S^3)$ and $X=\# n (T^2 \x S^2)$
for all $n > 0$. Also included are the compact, connected,
orientable 4-manifolds $X$ whose interiors $X - \bdry X$ admit a
complete hyperbolic metric. In addition, the corollary applies to
the total space of any orientable fiber bundle $S^2 \to X \to
\Sigma$ for some compact, connected, orientable 2-manifold $\Sigma$
of positive genus.
\end{exm}

\begin{rem}
Many surgical theorems on $\TOP$ 4-manifolds require $\pi$ to have
subexponential growth \cite{FreedmanTeichner, KrushkalQuinn} in
order to find topologically embedded Whitney discs.  Currently, the
Topological Surgery Conjecture remains open for the more general
class of discrete, amenable groups. In our case, observe that all
crystallographic groups $1 \to \Z^m \to \pi \to finite \to 1$ have
subexponential growth for all $m>0$. On the other hand, observe that
all torsion-free lattices $\pi$ in $\Isom(\qH^m)$ and all free
groups $\pi = F_n$ have exponential growth if and only if $m,n>1$.

Indeed, taking all $\qE^{m_i}=\R$, we obtain the finite-rank free
groups $\pi=F_n$.  Thus we partially strengthen a theorem of
Krushkal and Lee \cite{KrushkalLee} if $X$ is a compact, connected,
oriented $\TOP$ 4-manifold with fundamental group $F_n$. They only
required $X$ to be a finite Poincar\'e complex of dimension 4
($\bdry X = \varnothing$) but insisted that the intersection form
over $\Z[\pi]$ of their degree one, $\TOP$ normal maps $f: M \to X$
be tensored up from the simply-connected case $\Z[1]$.  Now, our
shortcoming is that exactness is not proven at $\cN_\TOP(X)$. This
is because self-equivalences do not represent the homotopy
equivalences with $\ks \neq 0$, such as the well-known
non-smoothable homotopy equivalences $*\CP^2 \to \CP^2$ and $*\RP^4
\to \RP^4$.
\end{rem}

Consider examples of infinite groups with odd torsion and trivial
orientation. The original case $n=1$ below was observed by S.
Cappell \cite[Thm. 5]{CappellSplit}. Observe that the free products
below have exponential growth if and only if $n>1$.

\begin{cor}\label{Cor_FreeOdd}
Let $\pi$ be a free product of groups of the form
\[
\pi = \bigstar_{i=1}^n O_i
\]
for some $n > 0$, where each $O_i$ is an odd-torsion group.
(Necessarily $\omega$ is trivial.) Then the surgery sequences
(\ref{SES_DIFF}) and (\ref{SES_TOP}) are exact.
\end{cor}

Consider non-orientable 4-manifolds $X$ whose fundamental group $\pi
= \bigstar n(\C_2)$ is infinite and has 2-torsion.  We denote $S^2
\rtimes \RP^2$ the total space of the 2-sphere bundle classified by
the unique homotopy class of non-nullhomotopic map $\RP^2 \to
BSO(3)$.  This total space was denoted as the sphere bundle
$S(\gamma \oplus \gamma \oplus \R)$ in the classification of
\cite{HKT}, where $\gamma$ is the canonical line bundle over
$\RP^2$.

\begin{cor}\label{Cor_OrderTwo}
Suppose $X$ is a $\DIFF$ 4-manifold of the form
\[
X = X_1 \# \cdots \# X_n \# r(S^2 \x S^2)
\]
for some $n > 0$ and $r \geq 0$, and each summand $X_i$ is either
$S^2 \x \RP^2$ or $S^2 \rtimes \RP^2$ or $\#_{S^1} n (\RP^4)$ for
some $1 \leq n \leq 4$. Then the surgery sequences (\ref{SES_DIFF})
and (\ref{SES_TOP}) are exact.
\end{cor}

In symplectic topology, the circle sum $M \#_{S^1} N$ is defined as
$(M - \overset{\circ}{E}) \cup_{\bdry E} (N - \overset{\circ}{E})$,
where $E$ is the total space of a 3-plane bundle over $S^1$ with
given embeddings in the 4-manifolds $M$ and $N$. The preceding
corollary takes circle sums along the order-two generator $\RP^1$ of
$\pi_1(P_j)$; the normal sphere bundle $\bdry E = S^2 \rtimes \RP^1$
is non-orientable. Observe that all the free products $\pi$ in
Corollary \ref{Cor_OrderTwo} have exponential growth if and only if
$n > 2$.

Next, consider non-orientable 4-manifolds $X$ whose fundamental
groups $\pi$ are infinite and torsion-free.  Interesting examples
have $\pi$ in Waldhausen's class $\Cl$ of groups with vanishing
Whitehead groups $\Wh_*(\pi)$ \cite[\S19]{WaldhausenKrings}, such as
Haken 3-manifold bundles over the circle.

\begin{cor}\label{Cor_Haken}
Suppose $X$ is a $\TOP$ 4-manifold of the form
\[
X = X_1 \# \cdots \# X_n \# r (S^2 \x S^2)
\]
for some $n>0$ and $r \geq 0$, and each summand $X_i$ is the total
space of a fiber bundle
\[
H_i \longra X_i \longra S^1.
\]
Here, we suppose $H_i$ is a compact, connected 3-manifold such that:
\begin{enumerate}
\item
$H_i$ is $S^3$ or $D^3$, or

\item
$H_i$ is irreducible with non-zero first Betti number.
\end{enumerate}
Moreover, if $H_i$ is non-orientable, we assume that the quotient
group $H_1(H_i;\Z)_{(\alpha_i)_*}$ of coinvariants is 2-torsionfree,
where $\alpha_i: H_i \to H_i$ is the monodromy homeomorphism. Then
the surgery sequence (\ref{SES_TOP}) is exact.
\end{cor}

Finally, consider certain surface bundles over surfaces, which have
fundamental group in the same class $\Cl$. Let $Kl = \RP^2 \# \RP^2$
be the Klein bottle, whose fundamental group $\pi^\omega = \Z^+
\rtimes \Z^-$ has the indicated orientation. Observe that any
non-orientable, compact surface of positive genus admits a collapse
map onto $Kl$.

\begin{cor}\label{Cor_Klein}
Suppose $X$ is a $\TOP$ 4-manifold of the form
\[
X = X_1 \# \cdots \# X_n \# r(S^2 \x S^2)
\]
for some $n > 0$ and $r \geq 0$, and each summand $X_i$ is the total
space of a fiber bundle
\[
\Sigma_i^f \longra X_i \longra \Sigma_i^b.
\]
Here, we suppose the fiber and base are compact, connected
2-manifolds, $\Sigma_i^f \neq \RP^2$, and $\Sigma_i^b$ has positive
genus. Moreover, if $X_i$ is non-orientable, we assume that the
fiber $\Sigma_i^f$ is orientable and that the monodromy action of
$\pi_1(\Sigma_i^b)$ of the base preserves any orientation on the
fiber (i.e. the bundle is simple). Then the surgery sequence
(\ref{SES_TOP}) is exact.
\end{cor}

\subsection{Proofs in the orientable case}

\begin{proof}[Proof of Theorem \ref{Thm_SES} for orientable $X$]
Suppose $X$ satisfies Hypothesis \ref{Hyp_Orient}. Let $f: M \to X$
be a degree one, normal map of compact, connected, oriented $\TOP$
4-manifolds such that: $\bdry f = \Id_{\bdry X}$ on the boundary,
$f$ has vanishing surgery obstruction $\sigma_*(f)=0$, and $f$ has
vanishing Kirby-Siebenmann stable $\PL$ triangulation obstruction
$\ks(f)=0$.

Then, by Proposition \ref{Prop_SelfEquivalence}, $f$ is $\TOP$
normally bordant to a homotopy self-equivalence $h: X \to X$
relative to $\bdry X$. Thus exactness is proven at $\cN_{\TOP
0}(X)$.

Note, since $X$ is orientable (\ref{Eqn_H0}), that
\[
\Tor_1(H_0(\pi;\Z^\omega),\Z_2) = \Tor_1(\Z,\Z_2) = 0.
\]
Then, by Proposition \ref{Prop_RedTOP}, $\red_\TOP$ induces an
isomorphism from $\cN_\DIFF(X)$ to $\cN_{\TOP 0}(X)$. Thus exactness
is proven at $\cN_\DIFF(X)$.
\end{proof}

\begin{proof}[Proof of Corollary \ref{Cor_Wall}]
The result follows immediately from Theorem \ref{Thm_SES}, since
\[
\kappa_2: H_2(\pi;\Z_2) =0 \longra L_4^h(\Z[\pi])
\]
is automatically injective.
\end{proof}

\begin{proof}[Proof of Corollary \ref{Cor_FreeProdCrys}]
By Theorem \ref{Thm_SES}, it suffices to show $\kappa_2$ is
injective by induction on $n$. Suppose $n=0$. Then it is
automatically injective:
\[
\kappa_2: H_2(1;\Z_2) = 0 \longra L_4^h(\Z[1]) = \Z.
\]

Let $\Lambda$ be a torsion-free lattice in either $\Isom(\qE^m)$ or
$\Isom(\qH^m)$ or $\Isom(\CH^m)$. Since isometric quotients of the
homogeneous\footnote{\textbf{Homogeneous:} the full group of
isometries acts transitively on the riemannian manifold.} spaces
$\qE^m$ or $\qH^m$ or $\CH^m$ have uniformly bounded curvature
matrix (hence $A$-regular), by \cite[Proposition
0.10]{FJ_TorsionFreeGL}, the connective (integral) assembly map is
split injective:
\[
A_\Lambda\!\gens{1}: H_4(\Lambda;\GTOP) = H_0(\Lambda;\Z) \oplus
H_2(\Lambda;\Z_2) \longra L_4^h(\Z[\Lambda]).
\]
The decomposition of the domain follows from the Atiyah-Hirzebruch
spectral sequence for the connective spectrum $\GTOP =
\qL.\!\gens{1}$. Therefore the integral lift of the 2-local
component is injective:
\[
\kappa_2: H_2(\Lambda;\Z_2) \longra L_4^h(\Z[\Lambda]).
\]

Suppose for some $n > 0$ that $\kappa_2$ is injective for $\pi_n =
\bigstar_{i=1}^{n-1} \Lambda_i$.  Let $\Lambda_n$ be a torsion-free
lattice in either $\Isom(\qE^{m_n})$ or $\Isom(\qH^{m_n})$. Write
$\pi := \pi_n * \Lambda_n$. By the Mayer-Vietoris sequence in
$K$-theory \cite{WaldhausenKrings}, and since
\[
\Wh_1(1)=\Wh_0(1)=0=\wt{\Nil}_0(\Z[1];\Z[\pi_n - 1],
\Z[\Lambda_n - 1]),
\]
note $\Wh_1(\pi) = \Wh_1(\pi_n) \oplus \Wh_1(\Lambda_n)$. Also,
since the trivial group $1$ is square-root closed in the
torsion-free groups $\pi$ and $\Lambda_n$, we have
\[
\UNil_4^h(\Z[1];\Z[\pi_n - 1], \Z[\Lambda_n - 1]) = 0
\]
by \cite[Corollary 4]{CappellUnitary}, which was proven in
\cite[Lemmas II.7,8,9]{CappellSplit}. So
\begin{eqnarray*}
L_4^h(\Z[\pi]) &=& \wt{L}_4^h(\Z[\pi_n]) \oplus L_4^h(\Z[1])
\oplus \wt{L}_4^h(\Z[\Lambda_n])\\
H_2(\pi;\Z_2) &=& H_2(\pi_n;\Z_2) \oplus H_2(1;\Z_2) \oplus
H_2(\Lambda_n;\Z_2)
\end{eqnarray*}
by the Mayer-Vietoris sequences in $L$-theory \cite[Theorem
5(ii)]{CappellUnitary} and group homology \cite[\S VII.9]{Brown}.
Therefore, since $\kappa_2$ factors through the summand
$\wt{L}_4^h(\Z[-])$, we conclude that
\[
\kappa_2 = \SmMatrix{\kappa_2 & 0 & 0\\ 0 & \kappa_2 & 0\\ 0 & 0 &
\kappa_2}: H_2(\pi;\Z_2) \longra L_4^h(\Z[\pi])
\]
is injective. The corollary is proven for $n$ factors, thus
completing the induction.
\end{proof}

\begin{proof}[Proof of Corollary \ref{Cor_FreeOdd}]
Since each $O_i$ is odd-torsion, a transfer argument \cite{Brown}
shows that $H_2(O_i;\Z_2)=0$. Then, by the Mayer-Vietoris sequence
in group homology \cite[\S VII.9]{Brown} and induction, we conclude
$H_2(\bigstar_{i=1}^n O_i;\Z_2)=0$. Therefore $\kappa_2$ is
automatically injective.
\end{proof}

\subsection{Proofs in the non-orientable case}

\begin{proof}[Proof of Theorem \ref{Thm_SES} for non-orientable $X$]
Suppose $X$ satisfies Hypothesis \ref{Hyp_NonorientFin}. Let $f: M
\to X$ be a degree one, $\TOP$ normal map such that $\sigma_*(f)=0$
and $\ks(f)=0$.  Then, by Proposition \ref{Prop_SelfEquivalence},
$f$ is $\TOP$ normally bordant to a homotopy self-equivalence $h: X
\to X$ relative to $\bdry X$.  Thus exactness is proven at
$\cN_{\TOP 0}(X)$.

Further suppose the non-orientable 4-manifold $X$ is smooth.  Since
in this case we assume that $\pi$ has an orientation-reversing
element of order two, by \cite[Theorem 3.1]{CS2}, there exists a
closed $\DIFF$ 4-manifold $X_1$ and a simple homotopy equivalence
$h_1: X_1 \to X$ such that $\eta_\DIFF(h_1) \neq 0$ and
$\eta_\TOP(h_1) = 0$. The above argument of Proof \ref{Thm_SES} in
the orientable case shows for non-orientable $X$ that the kernel of
$\cN_\DIFF(X) \to \cN_\TOP(X)$ is cyclic of order two. Note
\begin{eqnarray*}
\eta_\DIFF(h \circ h_1) &= \eta_\DIFF(h) + h_*\, \eta_\DIFF(h_1)
&\neq
\eta_\DIFF(h)\\
\eta_\TOP(h \circ h_1) &= \eta_\TOP(h) + h_*\, \eta_\TOP(h_1) &=
\eta_\TOP(h)
\end{eqnarray*}
by the surgery sum formula given in Proposition \ref{Prop_GPL} and
\cite[Proposition 4.3]{RanickiAlgebraicII}. Therefore $f$ is $\DIFF$
normally bordant to either the simple homotopy equivalence $h: X \to
X$ or $h \circ h_1: X_1 \to X$ relative to $\bdry X$. Thus exactness
is proven at $\cN_\DIFF(X)$.

Suppose $X$ satisfies Hypothesis \ref{Hyp_NonorientInf}. Let $f: M
\to X$ be a degree one, $\TOP$ normal map such that $\sigma_*(f)=0$
and $\ks(f)=0$.  Then, by Proposition \ref{Prop_SelfEquivalence},
$f$ is $\TOP$ normally bordant to a homotopy self-equivalence $h: X
\to X$ relative to $\bdry X$.  Thus exactness is proven at
$\cN_{\TOP 0}(X)$.
\end{proof}

The following formula generalizes an analogous result of J. Shaneson
\cite[Prop. 2.2]{Shaneson_Splitting5}, which was stated in the
smooth case.

\begin{prop}\label{Prop_GPL}
Suppose $M, N, X$ are compact $\PL$ manifolds. Let $f: M \to N$ be a
degree one, $\PL$ normal map such that $\bdry f: \bdry M \to \bdry
N$ is the identity map. Let $h: N \to X$ be a homotopy equivalence
such that $\bdry h: \bdry N \to \bdry X$ is the identity. Then there
is a sum formula for $\PL$ normal invariants:
\[
\wh{h \circ f} = \eta(h) + h_*(\wh{f}) \in [X/\bdry X, \G/\PL]_0.
\]
\end{prop}

\begin{proof}
Any element of the abelian group $[X/\bdry X, \G/\PL]_0$ is the
stable equivalence class of a pair $(\xi, t)$, where $\xi$ is a
$\PL$ fiber bundle over $X/\bdry X$ with fiber $(\R^n,0)$ for some
$n$, and $t: \xi \to \varepsilon^n = (X/\bdry X) \x (\R^n,0)$ is a
fiber homotopy equivalence of the absolute fiber $\R^n - \set{0}
\simeq S^{n-1}$. The abelian group structure on $[X/\bdry X,
\G/\PL]_0$ is the $\pi_0$ of the Whitney sum $H$-space structure on
the $\Delta$-set $\mathrm{Map}_0(X/\bdry X, G/\PL)$ defined
rigorously in \cite[Proposition 2.3]{RourkeHaupt}.

Let $\nu_M$ be the $\PL$ normal $(\R^n, 0)$-bundle of the unique
isotopy class of embedding $M \hookra S^{n+\dim(M)}$, where $n >
\dim(M)+1$. For a certain stable fiber homotopy trivialization $s$
induced by the embedding of $M$ and the normal map $f$, the normal
invariant of the degree one, $\PL$ normal map $(f,\xi): M \to N$ is
defined by
\[
\wh{(f,\xi)} = (\xi - \nu_N,s).
\]
For the homotopy equivalence $h: N \to X$ with homotopy inverse
$\ol{h}: X \to N$ and any $\PL$ bundle $\chi$ over $N$, define the
pushforward bundle $h_*(\chi) := (\ol{h})^*(\chi)$ over $X$. Let $r$
be the stable fiber homotopy trivialization associated to the degree
one, $\PL$ normal map $(h,h_*(\nu_N))$. Then note
\begin{eqnarray*}
\wh{h \circ f} &=& (h_*(\xi) - \nu_X,r+h_*(s))\\
&=& (h_*(\nu_N) - \nu_X,r) + (h_*(\xi) - h_*(\nu_N),h_*(s))\\
&=& \eta(h) + h_*(\wh{f}).
\end{eqnarray*}
Here, the addition is the Whitney sum of stable $\PL$ bundles with
fiber $(\R^n,0)$ equipped with stable fiber homotopy
trivializations.
\end{proof}

\begin{prop}[L\'opez de Medrano]\label{Prop_LdM}
The following map is an isomorphism:
\[
\kappa_2: H_2(\C_2;\Z_2) \longra L_4^h(\Z[\C_2]^-).
\]
\end{prop}

Note that the source and target of $\kappa_2$ are isomorphic to
$\Z_2$ \cite[Thm. 13A.1]{Wall}.

\begin{proof}
Observe that the connective assembly map
\[
A_\pi\!\gens{1}: H_\oplus(\pi;\GTOP^\omega) \longra
L_\oplus^h(\Z[\pi]^\omega)
\]
is a homomorphism of $L^*(\Z)$-modules. Then, by action of the
symmetric complex $\sigma^*(\CP^2) \in L^4(\Z)$, there is a
commutative diagram
\[
\begin{diagram}
\node{H_2(\C_2;\Z_2)} \arrow{e,t}{\kappa_2}
\arrow{s,tb}{(\Id;\,\xo\,\sigma^*(\CP^2))}{\iso}
\node{L_4^h(\Z[\C_2]^-)} \arrow{s,tb}{\xo\,\sigma^*(\CP^2)}{\iso}\\
\node{H_2(\C_2;\Z_2)} \arrow{e,t}{\kappa_2^{(8)}}
\node{L_8^h(\Z[\C_2]^-)}
\end{diagram}
\]
where the vertical maps are isomorphisms by decorated periodicity
\cite{Sullivan}.  So it is equivalent to show that $\kappa_2^{(8)}$
is non-trivial.

Consider the commutative diagram
\[
\begin{diagram}
\node{\cN_\PL(\RP^8)} \arrow{see,t}{\sigma_*} \arrow{s,t}{\red_\TOP}
\node{} \node{} \node{}
\node{H_2(\C_2;\Z_2)} \arrow{sww,t}{\kappa_2^{(8)}}\\
\node{\cN_\TOP(\RP^8)} \arrow{s,t}{\mathrm{transv}} \node{}
\node{L_8^h(\Z[\C_2]^-).} \node{}
\node{H_2(\RP^8;\pi_6(\GTOP)^-)} \arrow{n,b}{u_2}\\
\node{H^0(\RP^8;\GTOP)} \arrow[4]{e,t}{\cap [\RP^8]_{\qL^.}} \node{}
\node{} \node{} \node{H_8(\RP^8;\GTOP^-)} \arrow{n,b}{\proj}
\end{diagram}
\]
The $\PL$ surgery obstruction map $\sigma_*$ for $\RP^8$ was shown
to be non-trivial in \cite[Theorem IV.3.3]{LdM} and given by a
codimension two Kervaire-Arf invariant.  So the map $\kappa_2^{(8)}$
is non-trivial.  Therefore $\kappa_2$ is an isomorphism.
\end{proof}

\begin{proof}[Proof of Corollary \ref{Cor_OrderTwo}]
We proceed by induction on the number $n > 0$ of free $\C_2$ factors
in $\pi$ to show that $\kappa_2: H_2(\pi;\Z_2) \to
L_4^h(\Z[\pi]^\omega)$ is injective.

Suppose $n=1$.  Then
\[
\kappa_2: H_2(\C_2; \Z_2) \longra L_4^h(\Z[\C_2]^-)
\]
is an isomorphism by Proposition \ref{Prop_LdM}.

Suppose the inductive hypothesis is true for $n > 0$. Write
\begin{eqnarray*}
\pi_n &:=& \bigstar (n-1)(\C_2)\\
\pi^\omega &=&  (\pi_n)^{\omega_n} * (\C_2)^-.
\end{eqnarray*}
By the Mayer-Vietoris sequence in group homology \cite[\S
VII.9]{Brown}, we have
\[
H_2(\pi; \Z_2) = H_2(\pi_n;\Z_2) \oplus H_2(\C_2;\Z_2).
\]
By the Mayer-Vietoris sequence in $L^h_*$-theory \cite[Thm.
5(ii)]{CappellUnitary}, using the Mayer-Vietoris sequence in
$K$-theory \cite{WaldhausenKrings} for $h$-decorations, we have
\[
L_4^h(\Z[\pi]^\omega) = L_4^h(\Z[\pi_n]^{\omega_n}) \oplus L_4^h(\Z[\C_2]^-)
 \oplus \UNil_4^h(\Z; \Z[\pi_n - 1]^{\omega_n}, \Z^-).
\]
Since $\kappa_2$ is natural in groups with orientation character, we
have
\[
\kappa_2 = \SmMatrix{\kappa_2 & 0 & 0\\ 0 & \kappa_2 & 0} :
H_2(\pi;\Z_2) \longra L_4^h(\Z[\pi]^\omega).
\]
Therefore, by induction, we obtain that $\kappa_2$ is injective for
the free product $\pi^\omega$.

Let $i > 0$.  If $X_i = S^2 \x \RP^2$ or $X_i = S^2 \rtimes \RP^2$,
then a Leray-Serre spectral sequence argument shows that
\[
\Ker(u_2) = \Z_2 [S^2] \subseteq \Z_2 [S^2] \oplus \Z_2 [\RP^2] =
\Ker(v_2).
\]
If $X_i = P_j \#_{S^1} P_k$, then a Mayer-Vietoris and Poincar\'e
duality argument shows that
\[
\Ker(u_2) = \Z_2 ([\RP^2_j] + [\RP^2_k]) = \Ker(v_2).
\]
Hence Theorem \ref{Thm_SES} applies for both sets of $X_i$.
\end{proof}

\subsection{Proofs in both cases of orientability}

\begin{proof}[Proof of Corollary \ref{Cor_Haken}]
Write $\Lambda_i$ as the fundamental group and $\omega_i$ as the
orientation character of $X_i$.  Then the connective assembly map
\[
A_{\Lambda_i}\!\gens{1}: H_4(B\Lambda_i; \GTOP^{\omega_i}) \longra
L_4^h(\Z[\Lambda_i]^{\omega_i})
\]
is an isomorphism, as follows. Note $\Lambda_i = \pi_1(X_i) =
\pi_1(H_i) \rtimes \Z$.

Suppose $H_i$ has type (1). Then, by \cite[Theorem 13A.8]{Wall}, the
map $A_{\Lambda_i}\!\gens{1}$ is an isomorphism in dimension 4,
given by signature (mod 2 if $\omega_i \neq 1$).

Suppose $H_i$ has type (2). Then, by \cite[Theorem
1.1(1)]{Roushon_Haken} if $\bdry H_i$ is non-empty and by
\cite[Theorem 1.2]{Roushon_Betti} if $\bdry H_i$ is empty, the
connective assembly map $A_{\pi_1(H_i)}\!\gens{1}$ is an isomorphism
in dimensions 4 and 5. Since $\pi_1(H_i)$, hence $\Lambda_i$, is a
member of Waldhausen's class $\Cl$ \cite[Prop.
19.5(6,8)]{WaldhausenKrings}, we obtain $\Wh_*(\Lambda_i)=0$ by
\cite[Proposition 19.3]{WaldhausenKrings}. So, by the
Ranicki-Shaneson sequence in $L^h_*$-theory \cite[Thm.
5.2]{RanickiLIII}, and by the five-lemma, we obtain that the
connective assembly map $A_{\Lambda_i}\!\gens{1}$ is an isomorphism
in dimension 4.

Therefore, for both types, the integral lift $\kappa_2$ of the
2-local component of $A_{\Lambda_i}\!\gens{1}$ is injective.  So, by
the inductive Mayer-Vietoris argument of Corollary
\ref{Cor_FreeProdCrys}, we conclude that $\kappa_2$ is injective for
the free product $\pi = \bigstar_{i=1}^n \Lambda_i$.

If $X$ is orientable, then $X$ satisfies Hypothesis
\ref{Hyp_Orient}.  Otherwise, suppose $X$ is non-orientable. Then
consider all $X_i$ which are non-orientable.  If $H_i$ is
orientable, then the monodromy homeomorphism $\alpha_i: H_i \to H_i$
must reverse orientation.  So there is a lift $\pi_1(X_i) \to
\pi_1(S^1) \xra{\Id} \Z$ of the orientation character.  Otherwise,
if $H_i$ is non-orientable, then $H_1(X_i) = H_1(H_i)_{(\alpha_i)_*}
\x \Z$ by the Wang sequence and is 2-torsionfree by hypothesis.  So
there is a lift $\pi_1(X_i) \to H_1(X_i) \to \Z$ of the orientation
character. Hence there is an epimorphism $(\Lambda_i)^{\omega_i} \to
\Z^-$. Thus there is an epimorphism $\pi^\omega \to \Z^-$. So $X$
satisfies Hypothesis \ref{Hyp_NonorientInf}. Therefore, in both
cases of orientability of $X$, Theorem \ref{Thm_SES} is applicable.
\end{proof}

\begin{proof}[Proof of Corollary \ref{Cor_Klein}]
Write $\Lambda_i$ as the fundamental group and $\omega_i$ as the
orientation character of $X_i$.

Suppose $\Sigma_i^f = S^2$. Since $\pi_1(X_i) = \pi_1(\Sigma_i^b)$
is the fundamental group of an aspherical, compact surface, by the
proof of a result of J. Hillman \cite[Lemma
8]{Hillman_SurfaceBundles}, the connective assembly map
$A_{\Lambda_i}\!\gens{1}$ is an isomorphism in dimension 4.

Suppose $\Sigma_i^f \neq S^2$. Since $\Sigma_i^f$ and$\Sigma_i^b$
are aspherical, $X_i$ is aspherical. By a result of J. Hillman
\cite[Lemma 6]{Hillman_SurfaceBundles}, the connective assembly map
$A_{\Lambda_i}\!\gens{1}$ is an isomorphism in dimension 4.

Indeed, in both cases, the Mayer-Vietoris argument extends to fiber
bundles where the surfaces are aspherical, compact, and connected,
which are possibly non-orientable and with non-empty boundary (see
\cite[Thm. 2.4]{CHS} for detail).

Then the integral lift of the 2-local component of
$A_{\lambda_i}\!\gens{1}$ is injective:
\[
\kappa_2: H_2(\Lambda_i; \Z_2) \longra
L_4^h(\Z[\Lambda_i]^{\omega_i}).
\]
So, by the Mayer-Vietoris argument of Corollary
\ref{Cor_FreeProdCrys}, we conclude that $\kappa_2$ is an injective
for the free product $\pi = \bigstar_{i=1}^n \Lambda_i$.

If $X$ is orientable, then $X$ satisfies Hypothesis
\ref{Hyp_Orient}.  Otherwise, suppose $X$ is non-orientable. Then
consider all $X_i$ which are non-orientable. By hypothesis, the
fiber $\Sigma_i^f$ is orientable and the monodromy action of
$\pi_1(\Sigma_i^b)$ on $H_2(\Sigma_i^f;\Z)$, induced by the bundle
$\Sigma_i^f \to X_i \to \Sigma_i^b$, is trivial. We must have that
the surface $\Sigma_i^b$ is non-orientable.  So, since $\Sigma_i^b$
is the connected sum of a compact orientable surface and non-zero
copies of Klein bottles $Kl$, by collapsing to any $Kl$-summand,
there is a lift $\pi_1(X_i) \to \pi_1(\Sigma_i^b) \to \pi_1(Kl) \to
\Z$ of the orientation character. Hence there is an epimorphism
$(\Lambda_i)^{\omega_i} \to \Z^-$. Thus there is an epimorphism
$\pi^\omega \to \Z^-$. So $X$ satisfies Hypothesis
\ref{Hyp_NonorientInf}. Therefore, in both cases of orientability of
$X$, Theorem \ref{Thm_SES} is applicable.
\end{proof}

\subsection*{Acknowledgments}

The author is grateful to his doctoral advisor, Jim Davis, for
discussions on the assembly map in relation to smooth 4-manifolds
\cite{Davis}. The bulk of this paper is a certain portion of the
author's thesis \cite{Khan_Dissertation}, with various improvements
from conversations with Chris Connell, Ian Hambleton, Chuck
Livingston, Andrew Ranicki, John Ratcliffe, and Julius Shaneson.
Finally, the author would like to thank his pre-doctoral advisor,
Professor Louis H. Kauffman, for the years of encouragement and
geometric intuition instilled by him.

\bibliographystyle{alpha}
\bibliography{PostDissertation}

\end{document}